\documentclass[11pt,a4paper]{article}
\usepackage{amsmath,amssymb}
\usepackage{bm}
\usepackage{ascmac}
\usepackage[dvipdfmx]{graphicx}
\usepackage{enumerate}
\usepackage{tikz}
\usepackage[left=3cm, right=3cm]{geometry}
\usepackage{comment}

\newcounter{Cntr}
\usepackage{amsthm}
\newtheorem{thm}{Theorem}[section]

\newtheorem{lem}[thm]{Lemma}
\newtheorem{con}[thm]{Conjecture}

\newtheorem{prob}[thm]{Problem}

\newtheorem{cla}{Claim}[section]

\theoremstyle{definition}

\newtheorem*{prf}{Proof}
\newtheorem*{prf2}{Proof of Claims \ref{cla:added_1+} and \ref{cla:added_1}}

\usepackage{latexsym}
\def\qed{\hfill $\Box$}

\definecolor{mygray}{gray}{0.9}

\begin{document}

\title{Some conditions for hamiltonian cycles in 1-tough $(K_2 \cup kK_1)$-free graphs} 

\author{%
    Katsuhiro Ota\thanks{E-mail address: \texttt{ohta@math.keio.ac.jp}}
    %This work was supported by JSPS KAKENHI Grant Number 16H03952}
    \qquad Masahiro Sanka\thanks{E-mail address: \texttt{sankamasa@keio.jp}}\\ %This work was supported by JST Doctoral Program Student Support Project (JPMJSP2123).}\\
    Department of Mathematics, Keio University, \\
    3-14-1 Hiyoshi, Kohoku-ku, Yokohama 223-8522, Japan
} 

\date{}

\maketitle

\begin{abstract}
    Let $k \geq 2$ be an integer.
    We say that a graph $G$ is $(K_2 \cup kK_1)$-free if it does not contain $K_2 \cup kK_1$ as an induced subgraph.
    Recently, Shi and Shan conjectured that every $1$-tough and $2k$-connected $(K_2 \cup kK_1)$-free graph is hamiltonian.
    In this paper, we solve this conjecture by proving the statement; every 1-tough and $k$-connected $(K_2 \cup kK_1)$-free graph with minimum degree at least $\frac{3(k-1)}{2}$ is hamiltonian or the Petersen graph.
    %Also, we show that every $(k+1)$-connected $(K_2 \cup kK_1)$-free graph of order $n$ with minimum degree at least $\frac{3k-1}{2}$ and independence number strictly smaller than $\frac{n}{2}$ is hamiltonian-connected.
\end{abstract}


\section{Introduction}\label{sec:intro}

All graphs in this paper are finite, undirected and simple.
We use the terminology not defined here in \cite{Diestel-2017}, and additional definitions will be given as needed.
Let $G$ be a graph with vertex set $V(G)$ and edge set $E(G)$.
For $v \in V(G)$, we denote by $N_G(x)$ and $\deg_G(x)$ the neighborhood and the degree of $x$, respectively.
For $S \subset V(G)$, we define $N_G(S) = \bigcup_{x \in S} (N_G(x) \setminus S)$.
Also, for a subgraph $H$ of $G$, we write $N_G(H)$ for $N_G(V(H))$.
We denote $\delta(G)$, $\alpha(G)$ and $\omega(G)$ the minimum degree, the independence number and the number of components of $G$, respectively.
If $|V(G)|=1$, we say that $G$ is \emph{trivial}.
Let $R$ and $R'$ be two isomorphism classes of graphs and $k$ be a positive integer.
Then, we define $R \cup R'$ as the disjoint union of $R$ and $R'$, and $kR$ as the disjoint union of $k$ copies of $R$.
We use $K_n$ and $P_n$ to denote the complete graph and the path of order $n$, respectively.

We say that a graph $G$ is \emph{hamiltonian} if it contains a \emph{hamiltonian cycle}, that is, a spanning cycle of $G$.
The notion of toughness was introduced in the study of hamiltonian cycles by Chv\'{a}tal \cite{Chvatal-1973}.
The \emph{toughness} of a graph $G$, denoted by $t(G)$, is defined by
\[t(G) = \min \left\{\frac{|S|}{\omega(G-S)} \mid S \subset V(G) \text{ and } \omega(G-S) \geq 2 \right\},\]
or $t(G) = \infty$ if $G$ is complete.
For a real number $t \geq 0$, if $t(G) \geq t$, then we say that $G$ is $t$-\emph{tough}.
Clearly, a graph $G$ is $t$-tough if and only if $t \omega(G-S) \leq |S|$ for any subset $S \subset V(G)$ with $\omega(G-S) \geq 2$.
Also, it is easily shown that a non-complete $t$-tough graph $G$ is $\lceil 2t \rceil$-connected and $\alpha(G) \leq \frac{|V(G)|}{t+1}$.

It is well known that every hamiltonian graph is 1-tough, but the converse does not hold.
In 1973, Chv\'{a}tal \cite{Chvatal-1973} conjectured that the converse holds at least in an approximate sense.

\begin{con}\label{con:chvatal}
    There exists a constant $t_0 > 0$ such that every $t_0$-tough graph on at least three vertices is hamiltonian.
\end{con}

Bauer, Broersma and Veldman \cite{Bauer_et_al-2000} showed that for any $t < \frac{9}{4}$ there exists a non-hamiltonian $t$-tough graph.
So if Conjecture \ref{con:chvatal} is true, then such a constant $t_0$ must be at least $\frac{9}{4}$.
For details of known results on Conjecture \ref{con:chvatal}, we refer the reader to the survey \cite{Bauer_et_al-2006}.

For a graph $R$, we say that a graph $G$ is $R$-\emph{free} if it does not contain $R$ as an induced subgraph.
Several papers verified Conjecture \ref{con:chvatal} for $R$-free graphs for some graph $R$, for example, $2K_2$-free graphs \cite{Broersma_et_al-2014, Shan-2019, Ota_et_al-2022}, $(P_2 \cup P_3)$-free graphs \cite{Shan-2021}, $P_4$-free graphs \cite{Jung-1978, Li_et_al-2016}, $(P_3 \cup K_1)$-free graphs \cite{Nikog-2013, Li_et_al-2016}, and $(P_3 \cup 2K_1)$-free graphs \cite{Gao_et_al-2022}.

In this paper, we consider some conditions for $(K_2 \cup kK_1)$-free graphs to be hamiltonian, where $k \geq 2$ is an integer.
The known results so far for $(K_2 \cup kK_1)$-free graphs are listed below.

\begin{thm}[\cite{Li_et_al-2016}]\label{thm:Li}
    Every $1$-tough $(K_2 \cup 2K_1)$-free graph on at least three vertices is hamiltonian.
\end{thm}

\begin{thm}[\cite{Hatfield_et_al-2021}]\label{thm:Hat}
    Every $3$-tough $(K_2 \cup 3K_1)$-free graph on at least three vertices is hamiltonian.
\end{thm}

\begin{thm}[\cite{Shi_et_al-2022}]\label{thm:Shi}
    Let $k \geq 4$ be an integer and let $G$ be a $4$-tough and $2k$-connected $(K_2 \cup kK_1)$-free graph.
    Then $G$ is hamiltonian.
\end{thm}

Note that the toughness conditions in Theorems \ref{thm:Li}, \ref{thm:Hat} and \ref{thm:Shi} cannot be replaced by any connectivity conditions.
In fact, for any constant $m$, the complete bipartite graph $K_{m,m+1}$ is a non-hamiltonian $m$-connected $(K_2 \cup kK_1)$-free graph.
However, in \cite{Shi_et_al-2022}, Shi and Shan conjectured that the condition 4-tough of Theorem \ref{thm:Shi} is not sharp and can be relaxed to 1-tough.

\begin{con}[\cite{Shi_et_al-2022}]\label{con:K2kK1}
    Let $k \geq 4$ be an integer and let $G$ be a $1$-tough and $2k$-connected $(K_2 \cup kK_1)$-free graph.
    Then $G$ is hamiltonian.
\end{con}

In this paper, we solve Conjecture \ref{con:K2kK1} and generalize Theorems \ref{thm:Li} and \ref{thm:Hat} by proving the following stronger theorem.

\begin{thm}\label{thm:main1}
    Let $k \geq 2$ be an integer and let $G$ be a $k$-connected $(K_2 \cup kK_1)$-free graph with $\alpha(G) \leq \frac{|V(G)|}{2}$ and $\delta(G) \geq \frac{3(k-1)}{2}$.
    Then $G$ is hamiltonian or the Petersen graph.
\end{thm}

Note that for a $k$-connected $(K_2 \cup kK_1)$-free graph $G$, $\alpha(G) \leq \frac{|V(G)|}{2}$ if and only if $G$ is 1-tough.
To show the nontrivial implication of this statement, suppose that $G$ is not 1-tough.
Then, there exists $S \subset V(G)$ with
$\omega(G-S) \geq 2$ and $\omega(G-S) > |S|$.
Since $G$ is $k$-connected, we have $|S| \geq k$, and thus $\omega(G-S) \geq k+1$.
Therefore, since $G$ is $(K_2 \cup kK_1)$-free, each component of $G-S$ is trivial.
Thus, we find that $V(G) \setminus S$ is an independent set of order larger than $\frac{|V(G)|}{2}$.
Thus, Theorem \ref{thm:main1} implies that Conjecture \ref{con:K2kK1} is true, even if we replace ``$2k$-connected'' by ``$k$-connected and $\delta(G) \geq \frac{3(k-1)}{2}$'' with only one exception --- the Petersen graph.
Also, we note that when $k=2$, Theorem \ref{thm:main1} is equivalent to Theorem \ref{thm:Li}.

\begin{figure}[t]
    \centering
    \includegraphics{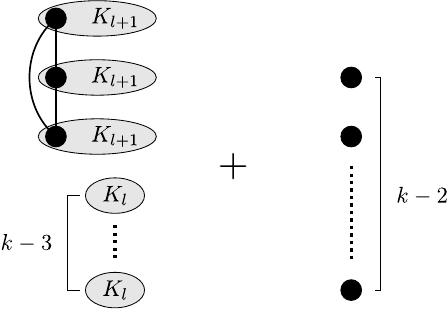}

    \caption{The graph $G(k,l)$}\label{fig:nonhamiltonian}
\end{figure}

When $k=3$, the connectivity condition of Theorem \ref{thm:main1} is sharp.
In fact, there exists a non-hamiltonian 1-tough (and thus, 2-connected) graph with independence number 3 and minimum degree arbitrarily large, which is shown in \cite{Li_et_al-2016}. 
Also, for each integer $k \geq 4$, there exists a non-hamiltonian $1$-tough and $(k-2)$-connected graph with independence number $k$ and minimum degree arbitrarily large.
To see this, for each pair of integers $k \geq 4$ and $l \geq 1$, we construct a graph $G(k,l)$ as follows.
For each $i \in \{1,2,3\}$, let $X_i$ be a copy of $K_{l+1}$, and $x_i \in V(X_i)$.
Let $X'=\left(\bigcup_{i=1}^3 X_i\right) + x_1 x_2 + x_2 x_3 + x_3 x_1$, and $X = X' \cup (k-3)K_l$.
Let $Y$ be a set of $k-2$ vertices with $X \cap Y = \emptyset$.
We define the graph $G(k,l)$, that is sketched in Figure \ref{fig:nonhamiltonian}, by $V(G(k,l)) = V(X) \cup Y$ and $E(G(k,l)) = E(X) \cup \{xy \mid x \in V(X) \text{ and }y \in Y\}$. 
Then the $G(k,l)$ is non-hamiltonian 1-tough and $(k-2)$-connected graph with $\alpha(G(k,l))=k$ and $\delta(G(k,l)) = k+l-3$.

For $l = 1$, the graph $G(k,1)$ is shown in \cite{Bigalke_et_al-1979}.
We note that the following theorem is proved in \cite{Bigalke_et_al-1979}.

\begin{thm}[\cite{Bigalke_et_al-1979}]
    Let $k \geq 4$ be an integer. 
    Then every $1$-tough and $(k-1)$-connected graph with independence number at most $k$ except for the Petersen graph is hamiltonian. 
\end{thm}

So we believe that when $k \geq 4$, the condition $k$-connected of Theorem \ref{thm:main1} can be relaxed to $(k-1)$-connected.
On the other hand, we are not sure whether the minimum degree condition of Theorem \ref{thm:main1} is sharp in general.
We are interested in whether the minimum degree condition can be relaxed to $k$.

\begin{prob}\label{newcon1}
    For each integer $k \geq 2$, is every $1$-tough and $k$-connected $(K_2 \cup kK_1)$-free graph except for the Petersen graph hamiltonian?
\end{prob}

If Problem \ref{newcon1} has a positive solution, then we obtain a generalization of the classical theorem of Chv\'{a}tal and Erd\H{o}s that for each $k \geq 2$, every $k$-connected graph with independence number at most $k$ is hamiltonian.

%%%%%%%%%%%%%%%%%

The rest of this paper is organized as follows.
In Section \ref{sec:pre}, we introduce some lemmas for the proof of Theorem \ref{thm:main1}. 
In Section \ref{sec:main1}, we prove Theorem \ref{thm:main1}.

\section{Preliminary}\label{sec:pre}

For the proof of Theorems \ref{thm:main1}, we introduce terminology and lemmas for (oriented) paths and cycles.

Let $P$ be a path.
For $x,y \in V(P)$, $xPy$ denotes the path from $x$ to $y$ passing through $P$.
%For a path $P$ with a given orientation, a vertex $x$ on $P$ and a positive integer $h$, $x^{+h}$ and $x^{-h}$ denote the $h$-th successor and the $h$-th predecessor of $x$ on $P$ (if they exist), respectively.
%We write $x^+$ and $x^-$ for $x^{+1}$ and $x^{-1}$, respectively.
Let $C$ be a cycle with a given orientation.
For $x \in V(C)$ and a positive integer $h$, $x^{+h}$ and $x^{-h}$ denote the $h$-th successor and the $h$-th predecessor, respectively.
We write $x^+$ and $x^-$ for $x^{+1}$ and $x^{-1}$, respectively.
For $u,v \in V(C)$, $u \overrightarrow{C} v$ denotes the $uv$-path from $u$ to $v$ along the orientation of $C$.
Also, $u \overleftarrow{C} v$ denotes $v \overrightarrow{C} u$.
For a subset $I \subset V(C)$, let $I^+ = \{x^+ \mid x \in I\}$ and $I^-=\{x^- \mid x \in I\}$.

We use the following lemma in the proof of Theorem \ref{thm:main1}.

\begin{lem}\label{lem:longest cycle}
    Let $G$ be a graph, $C$ be a longest cycle in $G$ with an orientation, and $H$ be a component of $G-C$.
    %Note that $N_G(H) \subset V(C)$.
    \begin{enumerate}[(i)]
        \item If $u \in N_G(H)$, then $u^+,u^- \notin N_G(H)$.
        \item Let $u,v \in N_G(H)$ be two vertices. 
        Then there are no $u^+v^+$-path and $u^-v^-$-path with all internal vertices in $V(G-C)$.
    \end{enumerate}
\end{lem}
\begin{prf}
    For (i), assume that $u^+ \in N_G(H)$.
    Then we can take a $uu^+$-path $P$ with all internal vertices on $H$.
    However, then
    $u P u^+ \overrightarrow{C} u$
    is a cycle of $G$ longer than $C$, a contradiction.
    Thus, we have $u^+ \notin N_G(H)$ and $u^- \notin N_G(H)$.

    For (ii), assume that there exists a $u^+v^+$-path $Q$ with all internal vertices in $V(G-C)$.
    Note that by (i), we have $V(Q) \cap V(H) = \emptyset$.
    Let $P$ be a $uv$-path with all internal vertices on $H$.
    Then
    $u P v \overleftarrow{C} u^+ Q v^+ \overrightarrow{C} u$
    is a cycle of $G$ longer than $C$, a contradiction.
    Thus, we find that (ii) holds. \qed
\end{prf}

In particular, every longest cycle of a $k$-connected $(K_2 \cup kK_1)$-free graph satisfies the following property.

\begin{lem}\label{lem:longest cycle of K2kK1}
    Let $k \geq 2$ be an integer and $G$ be a $k$-connected $(K_2 \cup kK_1)$-free graph.
    Then for every longest cycle $C$ of $G$, each component of $G-C$ is trivial.
\end{lem}
\begin{prf}
    We take a longest cycle $C$ of $G$ and give $C$ a fixed orientation.
    Let $H$ be a component of $G-C$.
    Then by Lemma \ref{lem:longest cycle}, the set $N_G(H)^+$ is independent, and no vertex in $H$ is adjacent to a vertex in $N_G(H)^+$.
    Also, since $G$ is $k$-connected, we have $|N_G(H)^+| = |N_G(H)|  \geq k$.
    Therefore, since $G$ is $(K_2 \cup kK_1)$-free, $H$ must be trivial. \qed
\end{prf}

Finally, we show the following key lemma for $(K_2 \cup kK_1)$-free graphs.

\begin{lem}\label{lem:K2kK1}
    Let $k \geq 2$ be an integer, $G$ be a $(K_2 \cup kK_1)$-free graph and $X \subset V(G)$ be an independent set.
    Then, the following statements hold.
    \begin{enumerate}[(i)]
        \item Each vertex $v \in V(G)$ satisfies either $N_G(v) \cap X = \emptyset$ or $|N_G(v) \cap X| \geq |X|-k+1$.
        \item Suppose that each vertex $w$ in $W \subset V(G)$ satisfies $|N_G(w) \cap X| \leq |X|-k$.
        Then, $X \cup W$ is independent in $G$. 
    \end{enumerate}
\end{lem}
\begin{prf}
    %For (1), take a vertex $v \in V(G)$ satisfying $N_G(v) \cap X \neq \emptyset$.
    %Since $G$ is $(K_2 \cup kK_1)$-free, we have $|X \setminus N_G(v)| \leq k-1$.
    %Therefore, we have $|N_G(v) \cap X| \geq |X|-k+1$.

    For (i), suppose to the contrary that $1 \leq |N_G(v) \cap X| \leq |X|-k$ for a vertex $v$.
    Let $x \in N_G(v) \cap X$ and $X' \subset X \setminus N_G(v)$ with $|X'| = k$.
    Then, $\{v,x\} \cup X'$ induces $K_2 \cup kK_1$, a contradiction.

    For (ii), suppose that $W \neq \emptyset$ and each vertex $w$ in $W \subset V(G)$ satisfies $|N_G(w) \cap X| \leq |X|-k$.
    Then, we have $|X| \geq k$.
    Also, by (i), there are no edges joining $W$ and $X$.
    Suppose that there exist $u,v \in W$ with $uv \in E(G)$.
    Then for a set $X' \subset X$ with $|X'|=k$, the set $\{u,v\} \cup X'$ induces $K_2 \cup kK_1$, a contradiction.
    Therefore, the set $X \cup W$ is independent.\qed 
\end{prf}

%%%%%%%%%%%%%%%%%%%%

\section{Proof of Theorem \ref{thm:main1}}\label{sec:main1}

Let $G$ be a $k$-connected $(K_2 \cup kK_1)$-free graph with $\delta(G) \geq \frac{3(k-1)}{2}$.
To prove Theorem \ref{thm:main1}, we show that if $G$ is neither hamiltonian nor the Petersen graph, then $\alpha(G)>\frac{|V(G)|}{2}$ holds.
Let $C$ be a longest cycle of $G$.
Since $G$ is not hamiltonian, we have $V(C) \neq V(G)$.
We give $C$ a fixed orientation.

Let $x \in V(G-C)$ and $d = \deg_G (x)$.
Since any component of $G-C$ is trivial by Lemma \ref{lem:longest cycle of K2kK1}, we have $d \geq \delta(G)$ and $N_G(x) \subset V(C)$.
Let $x_1,x_2, \ldots, x_d$ be the vertices in $N_G(x)$ appearing in this order along the orientation.
Let $X = \{x\} \cup N_G(x)^+ = \{x,x_1^+,\ldots,x_d^+\}$.
Then by Lemma \ref{lem:longest cycle}, $X$ is an independent set.
Also, we have $|X| = d + 1 \geq \delta(G)+1 \geq k$.
We show the following claims relating the longest cycle.

\begin{cla}\label{cla:2}
    Let $i,j \in \{1, \ldots,d\}$ with $i \neq j$ and let $uv \in E(x_j^+ \overrightarrow{C} x_i)$ with $v = u^+$.
    If $ux_i^+ \in E(G)$, then we have $vx_j^+ \notin E(G)$.
\end{cla}
\begin{prf}
    Suppose that $ux_i^+ \in E(G)$ and $vx_j^+ \in E(G)$.
    Then
    \[u x_i^+ \overrightarrow{C} x_j x x_i \overleftarrow{C} v x_j^+ \overrightarrow{C} u\]
    is a cycle in $G$ longer than $C$, a contradiction. \qed
\end{prf}

Next, we show the following claim relating the set $N_G(X)$.

\begin{cla}\label{cla:added_1+}
    For each edge $uv \in E(G)$, it follows that $|\{u,v\} \cap N_G(X)|=1$, unless $\{u,v\}=\{x_i^-, x_i\}$ for some $i \in \{1,\ldots,d\}$ and $\{u,v\} \subset N_G(X)$.
\end{cla}

The following claim describes the exceptional situation of Claim \ref{cla:added_1+}, where $\{x_i^-,x_i\} \subset N_G(X)$ for some $i$. 
We shall prove these claims simultaneously.

\begin{cla}\label{cla:added_1}
    Suppose that $x_i^- \in N_G(X)$ for some $i \in \{1,\ldots,d\}$.
    Then, $k$ is odd and $d= \frac{3(k-1)}{2}$.
    Moreover, the following conditions hold.
    \begin{enumerate}[(i)]
        \item $N_G(x_i^-) \cap X = \{x_{i+j}^+ \mid k-2 \leq j \leq d-1\}$;
        \item $N_G(x_i) \cap X = \{x\} \cup \{x_{i+j}^+ \mid 0 \leq j \leq d-k\}$; and
        \item $N_G(x^{+2}_{i+k-2}) \cap X = \{x_{i+j}^+ \mid d-k \leq j \leq k-2\}$,
    \end{enumerate}
    where the indices are taken modulo $d$.
\end{cla}
\begin{prf2}
    In order to prove these claims, by symmetry, we consider the edge in $x_d^+ \overrightarrow{C} x_1$.
    Let $uv \in E(x_d^+ \overrightarrow{C} x_1)$ with $v = u^+$.
    Note that $x_d^+ \in X \subset V(G) \setminus N_G(X)$ and $x_1 \in N_G(x) \subset N_G(X)$.
    Suppose $\{u,v\} \cap N_G(X) = \emptyset$.
    Then we have $uv \in E(x_{d}^{+3} \overrightarrow{C} x_1^-) \subset E(x_d^+ \overrightarrow{C} x_1)$ since $x_d^{+2} \in N_G(x_d^+) \subset N_G(X)$.
    However, for a set $X' \subset X$ with $|X'| = k$, the set $\{u,v\} \cup X$ induces $K_2 \cup kK_1$, a contradiction.
    Thus, we obtain $|\{u,v\} \cap N_G(X)| \geq 1$.
    
    Next, suppose $\{u,v\} \subset N_G(X)$.
    Then by Lemma \ref{lem:K2kK1} (i), we have
    \begin{equation}\label{eq1}
        \min\{|N_G(u) \cap X|, |N_G(v) \cap X|\} \geq |X|-k+1 = d - k + 2.
    \end{equation}
    Let $r \in \{1,\ldots,d\}$ be the maximum index such that $vx_r^+ \in E(G)$ and $l \in \{1,\ldots,d\}$ be the minimum index such that $u x_l^+ \in E(G)$. 
    Then by Claim \ref{cla:2}, we have $r \leq l$.
    Also, $N_G(u) \cap X \subset \{x_j^+ \mid l \leq j \leq d\}$ and $N_G(v) \cap X \subset \{x\} \cup \{x_j^+ \mid 1 \leq j \leq r\}$.
    These and the inequality (\ref{eq1}) imply 
    \begin{equation}\label{eq2}
    r \geq d-k+1 \text{ and } l \leq k-1.
    \end{equation}
    In particular, if $r = d-k+1$, then $v$ must be adjacent to $x$, and hence $v = x_1$ and $u = x_1^-$ hold.

    We consider the vertex $x_l^{+2}$.
    Since $x_l^{+2} \in N_G(x_l^+) \subset N_G(X)$, we have 
    \begin{equation}\label{eq3}
        |N_G(x_l^{+2}) \cap X| \geq d-k+2 \geq 2
    \end{equation}
    by Lemma \ref{lem:K2kK1} (i).
    If $x_l^{+2} = x_{l+1}$, then
    \[x x_{l+1} \overrightarrow{C} u x_l^+ \overleftarrow{C} x_r^+ v \overrightarrow{C} x_r x\]
    is a cycle in $G$ longer than $C$, a contradiction.
    Thus, we have $x_l^{+2} \neq x_{l+1}$, implying that there exists an index $j \in \{1,\ldots,d\} \setminus \{l\}$ such that $x_l^{+2} x_j^+ \in E(G)$.
    If $j < r$ or $j > l$, then
    \[  \left\{
        \begin{aligned}
            &x x_r \overleftarrow{C} x_j^+ x_l^{+2} \overrightarrow{C} u x_l^+ \overleftarrow{C} x_r^+  v \overrightarrow{C} x_j x & \text{ if } &j < r,\\
            &x x_j \overleftarrow{C} x_l^{+2} x_j^+ \overrightarrow{C} u x_l^+ \overleftarrow{C} x_r^+ v \overrightarrow{C} x_r x & \text{ if } &j > l
        \end{aligned}
        \right.
    \]
    is a cycle in $G$ longer than $C$, a contradiction.
    Therefore, we have $N_G(x_l^{+2}) \cap X \subset \{x_j^+ \mid r \leq j \leq l\}$, and thus $l-r+1 \geq d-k+2$ by (\ref{eq3}).

    Since $l-r+1 \leq k-1-(d-k) = 2k - d  - 1$ by (\ref{eq2}), we have $2d \leq 3(k - 1)$.
    Therefore, since $d \geq \delta(G) \geq \frac{3(k-1)}{2}$, we have the equality $d = \frac{3(k - 1)}{2}$, and hence $k$ is odd.
    Moreover, the inequalities in (\ref{eq2}) achieve the equalities, meaning $r = d-k+1$ and $l = k - 1$.
    Thus, we have $v = x_1$ and $u = x_1^-$.
    Also, we get $N_G(x_{k-1}^{+2}) \cap X = N_G(x_l^{+2}) \cap X = \{x_j^+ \mid r \leq j \leq l\} = \{x_{1+j}^+ \mid d-k \leq j \leq k-2\}$, $N_G(x_1^-) \cap X = \{x_{1+j}^+ \mid k-2 \leq j \leq d-1\}$ and $N_G(x_1) \cap X = \{x\} \cup \{x_{1+j}^+ \mid 0 \leq j \leq d-k\}$.
    These imply the conditions (iii), (i) and (ii), respectively. \qed 
\end{prf2}

Let $Y = \{x,x_1^-, x_2^-, \ldots, x_d^-\}$.
Then by Lemma \ref{lem:longest cycle}, $Y$ is an independent set.
Also, by applying Claims \ref{cla:added_1+} and \ref{cla:added_1} to the inverse orientation of $C$, we get the following claim.

\begin{cla}\label{cla:added_2}
    Let $i \in \{1,\ldots,d\}$.
    If $x_i^+ \in N_G(Y)$, then the following conditions hold.
    \begin{enumerate}[$(i)$]
        \item $N_G(x_i^+) \cap Y = \{x_{i+j}^- \mid 1 \leq j \leq d-k+2\}$;
        \item $N_G(x_i) \cap Y = \{x\} \cup \{x_{i+j}^- \mid k \leq j \leq d\}$; and
        \item $N_G(x^{-2}_{i+d-k+2}) \cap Y = \{x_{i+j}^- \mid d-k+2 \leq j \leq k\}$,
    \end{enumerate}
    where the indices are taken modulo $d$.
\end{cla}

Finally, by showing the following claim, we will confirm that $|V(C) \setminus N_G(X)| = \frac{|V(C)|}{2}$ holds.

\begin{cla}\label{cla:3}
    The vertices on $C$ are alternatingly contained in $V(G) \setminus N_G(X)$ and $N_G(X)$.
    Consequently, $|V(C)|$ is even.
\end{cla}
\begin{prf}
    If each subpath $x_{i-1}^+ \overrightarrow{C} x_i$ consists of even number of vertices, then the claim follows from Claim \ref{cla:added_1+}.
    Suppose that some subpath, say $x_d^+ \overrightarrow{C} x_1$, consists of odd number of vertices.
    Then, $x_1^- \in N_G(X)$, which implies by Claim \ref{cla:added_1} that $k$ is odd, $d = \frac{3(k-1)}{2}$, $N_G(x_1^-) \cap X = \{x_{j}^+ \mid k-1 \leq j \leq d\}$, $N_G(x_1) \cap X = \{x\} \cup \{x_{j}^+ \mid 1 \leq j \leq d-k+1\}$ and $N_G(x_{k-1}^{+2}) \cap X = \{x_{j}^+ \mid d-k+1 \leq j \leq k-1\}$.
    Thus, since $x_d^+ \in N_G(x_1^-) \subset N_G(Y)$, we have $N_G(x_d^+) \cap Y = \{x_j^- \mid 1 \leq j \leq d-k+2\}$ and $N_G(x_d) \cap Y = \{x\} \cup \{x_j^- \mid k \leq j \leq d\}$ by Claim \ref{cla:added_2} (i) and (ii).

    Suppose $k \geq 5$.
    Then we have $d > k$.
    Since $x_{k-1}^+ \in N_G(x_1^-) \subset N_G(Y)$, we get $N_G(x_{k-1}^+) \cap Y = \{x_{k-1+j}^- \mid 1 \leq j \leq d-k+2\}$ by Claim \ref{cla:added_2} (i).
    However, then
    \[x x_k \overrightarrow{C} x_d^- x_{k-1}^+ \overrightarrow{C} x_k^- x_d \overrightarrow{C} x_{k-1} x\]
    is a cycle in $G$ longer than $C$, a contradiction.
    Therefore, we have $k=3$, and hence $d=k=3$.

    Note that we have $N_G(x_1^-) \cap X= \{x_2^+,x_3^+\}$, $N_G(x_1) \cap X = \{x, x_1^+\}$, $N_G(x_2^{+2}) \cap X = \{x_1^+,x_2^+\}$, $N_G(x_3^+) \cap Y = \{x_1^-, x_2^-\}$ and $N_G(x_3) \cap Y = \{x, x_3^-\}$.
    Since $x_2^- \in N_G(X)$ and $x_2^+ \in N_G(Y)$, by Claims \ref{cla:added_1} (i) and \ref{cla:added_2} (i), we have $N_G(x_2^-) \cap X = \{x_1^+,x_3^+\}$ and $N_G(x_2^+) \cap Y= \{x_1^-,x_3^-\}$.
    Similarly, we get $N_G(x_1^+) \cap Y = \{x_2^-,x_3^-\}$ and $N_G(x_3^-) \cap X = \{x_1^+,x_2^+\}$ since $x_1^+ \in N_G(Y)$ and $x_3^- \in N_G(X)$, respectively.
    Note that these imply $x_1^+ \neq x_2^-$, $x_2^+ \neq x_3^-$ and $x_3^+ \neq x_1^-$.
    Also, by Claims \ref{cla:added_1} (iii) and \ref{cla:added_2} (iii), we have $N_G(x_1^{+2}) \cap X = \{x_1^+,x_3^+\}$ and $N_G(x_3^{-2}) \cap Y=\{x_1^-,x_3^-\}$.
    If $x_2^{+2} \neq x_3^-$, then
    \[x x_1 x_1^+ x_2^{+2} \overrightarrow{C} x_3^{-2} x_1^- \overrightarrow{C} x_3^+ x_1^{+2} \overrightarrow{C} x_2^+ x_3^- x_3 x\]
    is a cycle in $G$ longer than $C$, even if $x_1^{+2} = x_2^-$, a contradiction.
    Thus, we have $x_2^{+2} = x_3^-$.
    By the symmetry, we get $x_1^{+2}=x_2^-$ and $x_3^{+2}=x_1^-$.

    Note that by Claims \ref{cla:added_1} (ii) and \ref{cla:added_2} (ii), we have $N_G(x_i) \cap (X \cup Y)=\{x,x_i^-,x_i^+\}$ for each $i \in \{1,2,3\}$.
    Suppose that there exists a vertex $y \in V(G) \setminus (V(C) \cup \{x\})$.
    By Lemma \ref{lem:longest cycle} (ii), we have $\max\{|N_G(y) \cap X|, |N_G(y) \cap Y|\} \leq 1$.
    Thus, by Lemma \ref{lem:K2kK1} (i), $N_G(y)$ intersects neither $X$ nor $Y$.
    This and $\delta(G) \geq 3$ imply $N_G(y) = \{x_1,x_2,x_3\}$.
    However, if there exists an $x_1 x_3$-path $Q=x_1 y x_3$ or $Q=x_1 x_3$ in $G$, then
    \[x x_2 x_2^+ x_1^- x_3^+ x_2^- x_1^+ x_3^- x_3 Q x_1 x\]
    is a cycle in $G$ longer than $C$, a contradiction.
    Thus, we have $V(G) = V(C) \cup \{x\}$ and $x_1x_3 \notin E(G)$.
    Similarly, we get $x_1 x_2 \notin E(G)$ and $x_2 x_3 \notin E(G)$, and hence $N_G(x_i) = \{x,x_i^-,x_i^+\}$ for each $i \in \{1,2,3\}$.
    
    After all, we find that the vertex set and the edge set of $G$ are
    \[\begin{aligned}
        &V(G) = \{x\} \cup \{x_i^-,x_i,x_i^+ \mid i \in \{1,2,3\}\}, \text{ and }\\
        &E(G) = E(C) \cup \{xx_1, xx_2, xx_3, x_1^+x_3^-, x_2^+x_1^-, x_3^+x_2^- \}.
    \end{aligned}\]
    Thus, $G$ is the Petersen graph, a contradiction.
    This completes the proof of Claim \ref{cla:3}. \qed
\end{prf}

Let $W = V(C) \setminus N_G(X)$.
Then, by Claim \ref{cla:3}, we have $|W|=\frac{|V(C)|}{2}$.
Also, for each $y \in V(G-C) \cup W$, by Lemma \ref{lem:longest cycle} (ii), we have $|N_G(y) \cap X| \leq 1 \leq d-k+1 = |X|-k$.
Thus, by Lemma \ref{lem:K2kK1} (ii), the set $X \cup W \cup V(G-C) = W \cup V(G-C)$ is independent.
Therefore, we have
\[\alpha(G) \geq |W| + |V(G)|-|V(C)| = \frac{|V(C)|}{2} + |V(G)|-|V(C)|>\frac{|V(G)|}{2}.\]
This completes the proof of Theorem \ref{thm:main1}. \qed

%%%%%%%%%%%%%%%%%%%%%%%%%%%%%%%%

%%%%%%%%%%%%%%%%%

\section*{Acknowledgements}
%We would like to express our gratitude to the referees for their helpful comments.
The first author (Katsuhiro Ota) was supported by JSPS KAKENHI Grant Number 22KO3404.
The second author (Masahiro Sanka) was supported by JST Doctoral Program Student Support Project Grant Number JPMJSP2123.

%\cite{bauer_et_al-2000, Shan-2019, Shi_et_al-2022, Gao_et_al-2022, Bauer_et_al-2006}
%\cite{Li_et_al-2016, Zheng_et_al-2022, Shan-2021,  Ota_et_al-2022, Hatfield_et_al-2021, Diestel-2017, Chvatal-1973, Broersma_et_al-2014}

\bibliography{22K2kK1} %hoge.bibから拡張子を外した名前
\bibliographystyle{myplain} %hoge.bibから拡張子を外した名前

\end{document}